\documentclass[12pt,twoside]{article}
\usepackage{amsmath, amsthm, amsfonts}
\usepackage{amssymb}
\usepackage{amscd}
\usepackage{verbatim}

\begin{document}
\newcommand{\eich}{{\mathcal H}}
\newcommand{\HardyC}{{\left(\frac{m-2}{2}\right)^2}}
\newcommand{\betaq}{{nq/2-N}}
\newcommand{\betastar}{{2^*(1-m/2)}}
\newcommand{\TheSpace}{{{\cal D}^{1,2}}}
\newcommand{\cw}{\stackrel{D}{\rightharpoonup}}
\newcommand{\id}{\operatorname{id}}
\newcommand{\supp}{\operatorname{supp}}
\newcommand{\wlim}{\operatorname{w-lim}}
\newcommand{\mymu}{{|y|^{2^*(1-m/2)}}}
\newcommand{\R}{{\mathbb R}}
\newcommand{\N}{{\mathbb N}}
\newcommand{\Z}{{\mathbb Z}}
\newcommand{\Q}{{\mathbb Q}}
\newcommand{\la}{\label}
\newcommand{\xO}{\Omega}
\newcommand{\domain}{\R^{n+m}}
\newcommand{\Domain}{\R^{n+m}\setminus \R^{n}}
\newcommand\runninghead[1] {\pagestyle{myheadings}\markboth {{\footnotesize\it{\quad #1}\hfill}}{{\footnotesize\it{#1\hfill\quad}}}}\headsep=40pt
\newtheorem{theorem}{Theorem}[section]
\newtheorem{corollary}[theorem]{Corollary}
\newtheorem{lemma}[theorem]{Lemma}
\newtheorem{definition}[theorem]{Definition}
\newtheorem{remark}[theorem]{Remark}
\newtheorem{proposition}[theorem]{Proposition}
\newtheorem{conjecture}[theorem]{Conjecture}
\newcommand{\be}{\begin{equation}}
\newcommand{\ee}{\end{equation}}
\renewcommand\runninghead[2]{\pagestyle{myheadings}\markboth
{{\footnotesize\it{\quad #1}\hfill}}
{{\footnotesize\it{#2\hfill\quad}}}}\headsep=40pt

\runninghead{\hfill\today\hfill ---\hfill  A. Tertikas \& K.
Tintarev}{Minimizers for the Hardy-Sobolev-Maz'ya inequality}
\begin{titlepage}

\title{On existence of minimizers for the Hardy-Sobolev-Maz'ya
inequality}
\author{A.Tertikas\\
 {\small Department of Mathematics}\\ {\small University of Crete}\\
{\small 71409 Heraklion, Greece, and}\\
{\small Institute of  Applied and Computational Mathematics} \\
{\small FORTH, 71110 Heraklion, Greece }\\
{\small tertikas@math.uoc.gr}\\\and K.Tintarev
\\{\small Department of Mathematics}\\{\small Uppsala University}\\
{\small SE-751 06 Uppsala, Sweden}\\{\small tintarev@math.uu.se}}

\maketitle

\begin{abstract}
We show existence of minimizers for the Hardy-Sobolev-Maz'ya
inequality in $\R^{m+n}\setminus\R^n$ when either $m>2$, $n\ge 1$
or $m=1$, $n\ge 3$.
\end{abstract}

\vfill \noindent \small{The authors expresses their gratitude to
the faculties of mathematics departments at Technion - Haifa
Institute of Technology, at the University of Crete and of the University of Cyprus for their hospitality. A.T. acknowledges partial support by the RTN European
network Fronts--Singularities, HPRN-CT-2002-00274. K.T
acknowledges support as a Lady Davis Visiting Professor at
Technion and partial
support from University of Crete, University of Cyprus and Swedish Research Council.}

\end{titlepage}

\section{Introduction}

Let $N=n+m\ge 3$, $n=0,\dots,N-1$. The space $\R^N=\R^{m+n}$ will
be regarded here as a product $\R^n\times\R^m$ and the variables
in $\R^N$ will be denoted as $(x,y)$, $x\in \R^n$, $y\in\R^m$. In
this notations, the Hardy inequality involving the distance from
$\R^n\times\{0\}$ (which will be for brevity denoted as $\R^n$)
reads
\begin{equation}
\label{HS1} \int_{\domain}|\nabla u(x,y)|^2dxdy
\ge\HardyC\int_{\domain} \frac{u^2(x,y)}{|y|^{2}}dxdy, \quad u\in
C_0^\infty(\Domain).
\end{equation}
The constant $\HardyC$ that appears in (\ref{HS1}) is the best
constant and is not attained. Maz'ya [\cite{M}, Corollary 3, p.
97] was the first that discovered that an additional term with the
critical Sobolev exponent $2^*=\frac{2N}{N-2}$ can be added in the
right hand side. That is, when $n\neq 0$, the following
Hardy-Sobolev-Maz'ya inequality holds true:
\begin{equation}
\label{HS} \int_{\domain}|\nabla
u(x,y)|^2dxdy-\HardyC\int_{\domain}
\frac{u^2(x,y)}{|y|^{2}}dxdy\ge \kappa_{m,n}\|u\|^2_{2^*}, \quad
u\in C_0^\infty(\Domain).
\end{equation}

If $m=2$, the inequality (\ref{HS}) becomes the usual limit
exponent Sobolev inequality, since in this case
$C_0^\infty(\Domain)$  is dense in ${\cal D}^{1,2}(\domain)$. It
is worth noting that when $n=0$ and the distance is taken from the
origin, the inequality (\ref{HS1}) is no longer true (cf. Brezis
and V\'{a}zquez \cite{BV}). Let $\xO$ is a bounded domain with
$0\in\xO$, $X(r) = (1-\log r)^{-1}$ for $0<r\le 1$, and let
$D:=\sup_{x \in \xO}|x| < +\infty$. Then one has (\cite{FT},
Theorem A) the following analog of (\ref{HS}):
 \be
\la{15} \int_{\xO} |\nabla u|^2 dx -
    \left(\frac{N-2}{2} \right)^2    \int_{\xO} \frac{u^2}{|x|^2} dx  \geq
 C \left( \int_{\xO} |u|^{\frac{2N}{N-2}} X^{\frac{2(N-1)}{N-2}}\left(\frac{|x|}{D}\right) dx
 \right)^{\frac{N-2}{N}},
\qquad u \in C^{\infty}_0(\xO).
 \ee
Inequality (\ref{15}) involves the critical exponent, but contrary
to (\ref{HS}) it has a logarithmic correction. Moreover,
  it  is sharp in the sense that one
 cannot take a smaller power
of the logarithmic correction $X$.
 In this paper our interest is in the existence of minimizers to the Hardy-Sobolev-Maz'ya
 inequality (\ref{HS}).

In case of $m=1$ the set $\Domain$ is disconnected, so the problem
can be naturally restated as a problem on the half-space. However,
in order to keep the notations uniform, this reduction is not made
here. We exclude from consideration the case $m=2$ when the
inequality (\ref{HS}) becomes the usual Sobolev inequality with
the limit exponent and the case $m=N$ (that is, $n=0$) when the
inequality does not hold.

Note that the expression
\begin{equation}
\label{norm} \|u\|_{\eich_0}^2:=\int_{\domain}\left(|\nabla
u(x,y)|^2-\HardyC \frac{u^2(x,y)}{|y|^{2}}\right)dxdy
\end{equation}
is a quadratic form, positive definite due to (\ref{HS}) on
$C_0^\infty(\Domain)$, and therefore a scalar product. Also by
(\ref{HS}), the Hilbert space $\eich_0$, defined by completion of
$C_0^\infty(\Domain)$ with respect the norm above is continuously
imbedded into $L^{2^*}(\domain)$ whenever $n>0$, and the elements
of $\eich_0$ can be identified as measurable functions (modulo
a.e.).

Let $T:C_0^\infty(\Domain)\to C_0^\infty(\Domain)$ be given by \be
\label{transform} (Tv)(x,y)=|y|^{-\frac{m-2}{2}}v(x,y)\ee and
define a Hilbert space $\eich$ as a completion of
$C_0^\infty(\Domain)$ under the norm $\|v\|_\eich:=
\|Tv\|_{\eich_0}^2$. It is clear that $T$ extends to an isometry
between $\eich$ and $\eich_0$ as well as to an isometry between
$L^{2^*}(\domain,\mymu)$ and $L^{2^*}(\domain)$. In particular,
for $n\neq 0$ the space $\eich$ consists of measurable functions.
Furthermore, elementary computations show that

\be \label{normv} \|v\|_\eich:=
\|Tv\|_{\eich_0}^2=\int_{\domain}|y|^{-(m-2)}|\nabla v|^2dxdy, \ee

i.e. $\eich=\TheSpace(\domain,|y|^{-(m-2)})$, and  the inequality
(\ref{HS}) takes the equivalent form

\be \label{HSv}
 \int_{\domain} |y|^{-(m-2)}|\nabla v|^2dxdy\ge
\kappa_{m,n}\left(\int_{\domain}|y|^{\betastar}|v|^{2^*}dxdy\right)^{2/2^*}.
\ee


We prove the following statement:
\begin{theorem}
\label{main} The minimization problem
\begin{equation}
\label{minProb}\kappa_{m,n}=\inf_{\int_{\domain}|y|^{\betastar}|v|^{2^*}dxdy
=1}\int_{\domain} |y|^{-(m-2)}|\nabla v|^2dxdy
\end{equation}
has a point of minimum in $\eich=\TheSpace(\domain,|y|^{-(m-2)})$
whenever $m>2$, $n>0$ or $m=1$, $n\ge 3$.
\end{theorem}

\begin{theorem}
\label{main2} The minimization problem
\begin{equation}
\label{minProb2} \kappa_{m,n}=\inf_{u\in \eich_0:
\int_{\domain}|u|^{2^*}dxdy =1}\|u\|_{\eich_0}^2
\end{equation}
has a point of minimum in $\eich_0$, the completion of
$C_0^\infty(\Domain)$ with respect to the norm (\ref{norm}),
whenever $m>2$, $n>0$ or $m=1$, $n\ge 3$.
\end{theorem}

Due to the transformation (\ref{transform}) Theorems~\ref{main}
and \ref{main2} are equivalent. The minimizer of (\ref{minProb2})
resp. (\ref{minProb}) and the exact value of $\kappa_{m,n}$
remains unknown. Theorems~\ref{main} and \ref{main2} do not
include the case $m=1$ and $n=2$.

The problem (\ref{minProb}) is not compact, and we use a
concentration-compactness technique similar to one of
\cite{Brezis}, based on weak convergence argument and the
Brezis-Lieb lemma (\cite{BL}). Its application is, however, not
straightforward. The group of invariant transformations (which
include dilations) that suffices to treat a similar problem in
\cite{BadialeTarantello} or \cite{Tin-Halfspace}, does not suffice
here. In their case, the critical dilation invariance that is
caused by a singular weight, rather than by critical growth of
nonlinearity, reduces the nonlinear term to a subcritical one,
once the domain of the problem is partitioned into similar cells
(of varying diameter that goes to zero as the cell approaches the
singularity). However, reduction of the term $\int|u|^{2^*}$ to a
subcritical term requires a partition of both the domain and the
range of $u$ into similar compact sets, which makes it inevitable
to append the group of available invariant transformations by the
non-invariant translations in the $y$-variable. This is possible,
but only because the latter incur a variational penalty.

It might be useful for the reader more accustomed to the
P.-L.Lions' version of concentration compactness (\cite{PLL1a},
 \cite{PLL1b}, \cite{PLL2a}, \cite{PLL2b}) to give here some heuristic
interpretation of the problem in those terms. In the problem
(\ref{minProb}) four different types of concentration arise:
translations in the $x$-variable, translations in the
$y$-variable, concentration in the interior and concentration at
the boundary (including concentration at infinity). Translations
in $|y_k|\to\infty$ incurs a variational penalty and so does the
interior concentration, provided that the infimum value
$\kappa_{m,n}$ is less than the Sobolev constant

\be \label{SN} S_N:=\inf_{w\in {\mathcal D}^{1,2}(\R^N):
\|w\|_{2^*}=1}\int_{\R^N}|\nabla w|^2,\; N\ge 3.\ee

We have established that $\kappa_{m,n}<S_{m+n}$ whenever $m+n>3$.
The remaining concentrations, the concentration on the boundary
and the translations in $x$, are due to invariant transformations
and are handled by the subadditivity argument. \vskip5mm

In order to consider the analog of the  problem (\ref{HS}) on an
open set $\Omega\subset\Domain$ we would like to start with a
well-known Brezis-Nirenberg problem \cite{BN}. Set first
$\label{SNO} S_N(\Omega):=\inf_{w\in C_0^\infty(\Omega):
\|w\|_{2^*}=1}\int_{\Omega}|\nabla w|^2$, $N\ge 3$. It is well
known that for every $\Omega$, $S_N(\Omega)=S_N$ and that there is
no minimizer when $\Omega\neq\R^N$.

In \cite{BN} one considers a bounded set $\Omega\subset\R^N$ and
the minimization problem

\be S_{\lambda,N}(\Omega):=\inf_{u\in H_0^1(\Omega):
\|w\|_{2^*}=1}\int_{\Omega}(|\nabla u|^2-\lambda u^2),\ee where
$\lambda>0$ does not exceed the first eigenvalue $\lambda_1$ of
the Dirichlet Laplacian in $\Omega$. It is shown in \cite{BN} that
the inequality $S_{\lambda,N}(\Omega)<S_N$ (separation of the
infimum from the concentration level) holds for $N>3$ ( as well as
for $\lambda$ sufficiently close to $\lambda_1$ when $N=3$), from
which existence of the minimizer for $S_{\lambda,N}(\Omega)$
easily follows.

In our case we consider, for an open set
$\Omega\varsubsetneq\Domain$, $n\neq 0$, $m\neq 2$, $m+n\ge 3$,
the minimization problem

\begin{equation}
\label{minProb4a}\kappa_{m,n}(\Omega)=\inf_{u\in
C_0^\infty(\Omega):\;\int_{\Omega}|u|^{2^*}dxdy =1}\int_{\Omega}
\left(|\nabla u|^2-\HardyC\frac{u^2}{|y|^2}\right) dxdy.
\end{equation}

An equivalent problem under transformation (\ref{transform}) is

\begin{equation}
\label{minProb5a} \kappa_{m,n}(\Omega)=\inf_{u\in
\TheSpace(\Omega,|y|^{-(m-2)}):\;
\int_{\Omega}|y|^{\betastar}|v|^{2^*}dxdy =1}\int_{\Omega}
|y|^{-(m-2)}|\nabla v|^2dxdy,
\end{equation}
where $\TheSpace(\Omega;|y|^{-(m-2)})$ is the completion of
$C_0^\infty(\Omega)$ with respect to the norm $\left(\int_{\Omega}
|y|^{-(m-2)}|\nabla v|^2dxdy\right)^\frac12$.

There are still four types of concentration as in the case of
$\Domain$. Concentration in the interior yields
$S_N>\kappa_{m,n}(\Omega)$, provided that $N>3$. Concentration at
a boundary point with $y=0$ occurs at the energy level
$\kappa_{m,n}$. By monotonicity, for any $\Omega$,
$\kappa_{m,n}(\Omega)\ge\kappa_{m,n}$.

It is easy to see that whenever
$\kappa_{m,n}(\Omega)=\kappa_{m,n}$, $\Omega\neq\Domain$, the
constant $\kappa_{m,n}(\Omega)$ is not attained: if $v$ were a
minimizer for $\kappa_{m,n}(\Omega)$, it would be then a minimizer
for $\kappa_{m,n}$, contrary to the maximum principle. In
particular,

\begin{theorem}(Non-existence of minimizers for $\Omega\varsubsetneq\Domain$.)
\label{main6} Assume that $m\neq 2$, $n\neq 0$ and $m+n\ge 3$.
Then $\kappa_{m,n}(\Omega)=\kappa_{m,n}$  provided that one of the
following conditions holds true:
\begin{description}
\item[a)] there exist $x_0\in\R^n$ and $r>0$ such that $\Omega$
contains $B_r(x_0,0)\setminus\{y=0\}$, \item[b)] there exists
$R>0$ such that $\Omega$ contains the set $\{|y|>R\}$.
\end{description}
\end{theorem}

There are also domains where we have existence of a minimizer
(consequently $\kappa_{m,n}(\Omega)>\kappa_{m,n}$). We consider
existence only for domains $\Omega$ contained in \be A_{r,R}=
\{(x,y)\in\Domain,\; r<|y|<R\},\; 0<r<R<\infty.\ee This condition
is not a heavy restriction in view of Theorem~\ref{main6}. Under
this assumption there is no concentration related to translations
in $y$. Concentration due to translations in the $x$-variable is
handled by subadditivity, under a flask-type assumption on
$\Omega$, that is: for every sequence $x_k\in\R^n$, there exists
$x_0\in\R^n$ such that
\be \label{flask-x} \liminf(\Omega+x_k)\subset \Omega+x_0.\ee

The existence proof for non-invariant domains cannot use the
Brezis-Lieb lemma directly. Instead, following the method of
\cite{acc}, it uses Lemma~\ref{BLsum} - an "iterated" version of
Brezis-Lieb lemma.

\begin{theorem}
\label{main4} Suppose that for some $0<r<R<\infty$,
\be\label{anncyl} \Omega\subset A_{r,R},\ee $\partial\Omega\in
C^1$ and, in addition, $\Omega$ satisfies (\ref{flask-x}). Then
the minimization problem

\begin{equation}
\label{minProb4}\kappa_{m,n}(\Omega)=\inf_{\int_{\Omega}|y|^{\betastar}|v|^{2^*}dxdy
=1}\int_{\Omega} |y|^{-(m-2)}|\nabla v|^2dxdy
\end{equation}
attains a minimum in $H^1_0(\Omega)$ provided that $m>2$, $n>0$,
or $m=1$, $n\ge 3$.
\end{theorem}

\begin{remark} Under hypothesis (\ref{anncyl}) we have
$\TheSpace(\Omega,|y|^{-(m-2)})=H^1_0(\Omega)$. This in particular
implies existence of the minimizer of (\ref{minProb4a}) in the
class $H^1_0(\Omega)$. This is not the case when $\Omega=\Domain$.
\end{remark}

In Section 2 we make preliminary computations used later in the
proofs. In Section 3 we prove that minimization sequences under
unbounded translations in the $y$-variable converge weakly to
zero. Section 4 concludes the proof of the main result (Theorems
\ref{main} and \ref{main2}). In Section 5 we prove
Theorems~\ref{main6} and \ref{main4}. In Section 6 we give an alternative proof of the main result and outline some open problems.

In what follows, integration without domains or variables
specified will always refer to $\Domain$ and $dxdy$, respectively.

\section{Preliminary computations}

\begin{lemma}
\label{w-eps} Let $w\in H^1_{loc}(\Domain)$. For every
$\epsilon\in(1,\frac14)$ there exists a $w_\epsilon\in
C_0^\infty(\Domain)$ such that

\be \label{w-epsilon} \int\nabla w\cdot\nabla w_\epsilon \ge
(1-\epsilon)\int|\nabla w_\epsilon|^2.
 \ee
Moreover, if $w\in {\mathcal D}^{1,2}(\Domain)$ then $w_\epsilon$
satisfies, in addition to (\ref{w-epsilon}),
$\|w-w_\epsilon\|_{{\mathcal D}^{1,2}}\le\epsilon \|w\|_{{\mathcal
D}^{1,2}}$.
\end{lemma}

\begin{proof} Assume first that $w\in {\mathcal
D}^{1,2}(\Domain)$ and let $\epsilon>0$. In this proof we use the
notation of the norm and of the scalar product in reference to the
space ${\mathcal D}^{1,2}(\Domain)$. By density of
$C_0^\infty(\Domain)$ in ${\mathcal D}^{1,2}(\Domain)$, and since
$w\neq 0$, one can choose a $w_\epsilon\in C_0^\infty(\Domain)$
such that $\|w-w_\epsilon\|\le\epsilon\|w\|$.  Using the Cauchy
inequality, we have \be (w,w_\epsilon)=
\|w\|^2-(w,w-w_\epsilon)\ge \|w\|^2-\|w\|\|w-w_\epsilon\|\ge
(1-\epsilon)\|w\|^2.\ee This proves the second assertion of the
lemma.

It remains now to consider the case $w\in
H^1_{loc}(\Domain)\setminus{\mathcal D}^{1,2}(\Domain)
\setminus\{0\}$. Then

\be \sup_{\psi\in
C_0^\infty(\Omega):\,\|\psi\|=1}(w,\psi)=+\infty,\ee since the
finite value of the supremum yields $w\in {\mathcal
D}^{1,2}(\Domain)$. In particular, there exists a $w_1\in
C_0^\infty(\Domain)$, $\|w_1\|=1$, such that $(w,w_1)>1$. Set
$w_\epsilon=w_1$.
\end{proof}

\begin{lemma}
\label{kappa-S} Let $\Omega\subset\Domain$ be an open set. If

\begin{description}
    \item[(i)] $m>2$ and $n\ge 1$, or
    \item[(ii)] $m=1$ and $n\ge 3$.
\end{description}
then $0<\kappa_{m,n}< S_{m+n}$.
\end{lemma}
Trivially, $\kappa_{2,n}=S_{2+n}$. We do not know whether, in the
remaining case, $\kappa_{2,1}<S_3$ or the equality prevails.

\begin{proof} Let $z=(x,y)\in\R^N$. The unique minimizer for (\ref{SN}), modulo translations
and the scale transformation $w\mapsto
\epsilon^{-\frac{N-2}{2}}w(z/\epsilon)$ is the well known
Bliss-Talenti solution, a scalar multiple of
$w=(1+|z|^2)^{-\frac{N-2}{2}}$.

Case (i). When $m>1$, $C_0^\infty(\Domain)$ is dense in ${\mathcal
D}^{1,2}(\domain)$. Then, since $\|u\|_{\eich_0}\le
\|u\|_{\mathcal D}^{1,2}$, the space ${\mathcal D}^{1,2}(\domain)$
is continuously imbedded into $\eich_0$ and for every $u\in
{\mathcal D}^{1,2}(\domain)$, $\int \frac{u^2}{|y|^2}<\infty$ and
\be \|u\|_{\eich_0}^2=\int|\nabla u|^2-\HardyC\int
\frac{u^2}{|y|^2}.\ee Substitution of $u=w$ proves therefore that
$\kappa_{m,n}<
 S_N$ whenever $m>2$.

Case (ii). Let $z_0=(x_0,y_0)\in\Domain$ and let
$\rho\in(0,\frac{|y|}{3})$. Let $\psi\in
C_0^\infty(B_\rho(z_0);[0,1])$ be equal $1$ on $B_\rho(z_0)$.
 These parameters will remain fixed. Let now
$w_\epsilon=\epsilon^{-\frac{N-2}{2}}w((x-x+0)/\epsilon(y-y_0)/\epsilon)$.
It suffices to prove that for $\epsilon$ sufficiently small,

\be\dfrac{ \int|\nabla(\psi w_\epsilon)|^2 -\frac14\int\frac{(\psi
w_\epsilon)^2}{y^2} } {\left(\int (\psi
w_\epsilon)^{2^*}\right)^\frac{2}{2^*}}< S_N,
 \ee
since  the left hand side is greater or equal to $\kappa_{1,n}$.

Note that $y$ is bounded from above and from below on
$B_\rho(z_0)$, so it suffices to show that for every $\lambda>0$

\be \label{BN} \dfrac{\int_{B_\rho(z_0)}|\nabla(\psi
w_\epsilon)|^2 -\lambda\int_{B_\rho(z_0)}(\psi w_\epsilon)^2}
{\left(\int_{B_\rho(z_0)} (\psi
w_\epsilon)^{2^*}\right)^\frac{2}{2^*}}< \dfrac{\int_{\R^N}|\nabla
w_\epsilon|^2} {\left(\int_{\R^N}
w_\epsilon^{2^*}\right)^\frac{2}{p}}=S_N.
 \ee

Verification of this is a literal repetition of the argument in
(\cite{BN}), cases $N=4$ and $N>4$, and can be omitted.
\end{proof}

\begin{remark}
\label{kappa-S1} Let $\Omega\subset\Domain$ be an open set. Set
\be \kappa_{m,n}(\Omega):=\inf_{u\in C_0^\infty(\Omega):
\int_{\domain}|u|^{2^*} =1}\|u\|_{\eich_0}^2.\ee Then
$\kappa_{m,n}(\Omega)<S_N$ whenever $m,n$ as in
Lemma~\ref{kappa-S}. The proof follows literally that of
Lemma~\ref{kappa-S1}, part (ii), provided that the point $z_0$ is
chosen in $\Omega$. Note that for $m>2$, $n\neq 0$ one has always
$m+n>3$.
\end{remark}

\begin{definition} \label{def:cwconvergence} \rm
Let $H$ be a Hilbert space equipped with a group $G$ of bounded
operators. We say that a sequence $u_{k}\in X$ converges to $u\in
X$ $G$-weakly, which we will denote as
\[
u_{k}  \stackrel{G}{\rightharpoonup} u,
\]
if for every sequence $g_k\in G$,
\begin{equation}
\label{eq:cwconvergence} g_k(u_k-u)\rightharpoonup 0.
\end{equation}
\end{definition}

Consider the following group acting on $\Domain$:

\begin{equation}
\label{d} d:=\{ \eta_{\alpha,j}: (x,y)\mapsto (2^{-j}
x-\alpha,2^{-j}y),\; j\in\R, \alpha\in\R^{n} \}.
\end{equation}

We associate with the group $d$ the following group of unitary
operators on $\eich_0$:

\begin{equation}
\label{setD0} D_0:=\{g_{\alpha,j}: v \mapsto
2^{-j(N-2)/2}v\circ\eta_{\alpha,j}, \eta_{\alpha,j}\in d\}.
\end{equation}

Operators in $D_0$ also preserve the $L^{2^*}$-norm.

By the isometry (\ref{transform}) $D=T^{-1}D_0T$ defines a group
of unitary operators on $\eich$ (which also preserve $\int
|y|^{\betastar}|v|^{2^*}$):
\begin{equation}
\label{setD} D:=\{g_{\alpha,j}: v \mapsto
2^{-jn/2}v\circ\eta_{\alpha,j}, \eta_{\alpha,j}\in d\}.
\end{equation}

\section{Penalty at infinity}

\begin{lemma}
\label{lem:localbound} Let $u_k$ be a bounded sequence in
$\eich_0$. If $|y_k|\to\infty$, then for all $k$ sufficiently
large $u_k(\cdot+(0,y_k))$ is bounded in $H^1_{loc}(\domain)$.
\end{lemma}
\begin{proof} Let $\Omega\subset\domain$ be an open bounded
set. Then by the H\"older inequality

\be \int_\Omega |y-y_k|^{-2}|u_k(\cdot+(0,y_k))|^2\le
\left(\int_\Omega |u_k(\cdot+(0,y_k))|^{2^*}\right)^\frac{2}{2^*}
\left(\int_\Omega |y-y_k|^{-N}\right)^\frac{2}{N}. \ee

The first integral in the right hand side is bounded since $u_k$
is bounded in $\eich_0$ and, therefore, by (\ref{HS}) in
$L^{2^*}$. The expression under the second integral converges
uniformly to zero. Therefore, the left hand side converges to
zero, and consequently,

\begin{eqnarray*}
C\ge\int_\Omega |y+y_k|^{2-m}|\nabla (|y+y_k|^{\frac{m-2}{2}}
u_k(\cdot+(0,y_k))|^2\ge\\\frac12 \int_\Omega|\nabla
u_k(\cdot+(0,y_k))|^2-C\int_\Omega
|y+y_k|^{-2}u_k(\cdot+(0,y_k))^2=\\ \int_\Omega\frac12 |\nabla
u_k(\cdot+(0,y_k))|^2+o(1).
\end{eqnarray*}

Therefore $\int_\Omega\frac12 |\nabla u_k(\cdot+(0,y_k))|^2$ is
bounded. It remains to note that
$\|u_k(\cdot+(0,y_k)\|_{2^*}=\|u_k\|_{2^*}$, which is bounded by
the $\eich_0$-norm.
\end{proof}

We call the sequence $u_k\in\eich_0$ (resp. $v_k\in\eich$) a
minimizing sequence, if $\|u_k\|_{2^*}=1$ and
$\|u_k\|_{\eich_0}^2\to\kappa_{m,n}$ (resp.
$\|v_k\|_{2^*,\mymu}=1$ and $\|v_k\|_{\eich}^2\to\kappa_{m,n}$).

\begin{lemma}
\label{runaway} Assume that $\kappa_{m,n}<S_N$. If $u_k\in\eich_0$
is a minimizing sequence and $|y_k|\to\infty$, then
$u_k(\cdot+(0,y_k))\rightharpoonup 0$ in $H^1_{loc}(\Domain)$ and
in $L^{2^*}(\domain)$.
\end{lemma}
\begin{proof} If the assertion of the lemma
is false, then there is a $w\in L^{2^*}\setminus\{0\}$ and (taking
into account Lemma~\ref{lem:localbound}) a renumbered subsequence
such that $u_k(\cdot+(0,y_k))\rightharpoonup w$ in $H^1_{loc}$ and
in $L^{2^*}$. Assume now that, on a renumbered subsequence, $\int
|u_k(\cdot+(0,y_k))-w|^{2^*}\to t\in[0,1]$.

Assume that $t\neq 1$. By the Brezis-Lieb lemma for
$L^{2^*}(\domain)$ (\cite{BL}),

\be \label{bll} \int_{\domain} |w|^{2^*}=1-t\ee Let $w_\epsilon\in
C_0^\infty(\Domain)$ be given by Lemma~\ref{w-eps} and let
$v^\epsilon_k:=u_k-w_\epsilon(\cdot-(0,y_k))$. Observing that,
since $w_\epsilon$ has compact support, \be
\left|\int|y|^{-2}u_kw_\epsilon(\cdot-(0,y_k))\right|=
\left|\int|y+y_k|^{-2}u_k w_\epsilon\right|\le
C_\epsilon|y_k|^{-2}\|u_k\|_{\eich_0}\to 0,
 \ee

and

\be \int|y|^{-2}w_\epsilon(\cdot-(0,y_k))^2=
\int|y+y_k|^{-2}{w_\epsilon}^2\to0, \ee

we have the following estimate:

\begin{eqnarray*}
\kappa_{m,n}=\|u_k\|_{\eich_0}^2+o(1)=&&\\
\|v_k^\epsilon\|_{\eich_0}^2+ \int|\nabla {w_\epsilon}|^2+
2\int\nabla v_k^\epsilon\cdot \nabla
w_\epsilon(\cdot-(0,y_k))-2\HardyC \int |y|^{-2}v_k^\epsilon
w_\epsilon(\cdot-(0,y_k)))+o(1)=&&\\
\|v_k^\epsilon\|_{\eich_0}^2+\int|\nabla w_\epsilon|^2+ 2\int
\nabla u_k\cdot \nabla w_\epsilon(\cdot-(0,y_k))+&&\\
2\HardyC\int|y|^{-2}u_kw_\epsilon(\cdot-(0,y_k)) -2\int|\nabla
w_\epsilon|^2-2\HardyC\int |y+y_k|^{-2}|w_\epsilon|^2+o(1)=&&\\
\|v_k^\epsilon\|_{\eich_0}^2+\int|\nabla w_\epsilon|^2+
2\int\nabla w\cdot\nabla w_\epsilon-\int|\nabla
w_\epsilon|^2+o(1)\ge
\|v_k^\epsilon\|_{\eich_0}^2+(1-2\epsilon)\int|\nabla
w_\epsilon|^2+o(1).&&
\end{eqnarray*}

Note that from this estimate follows that that $\int|\nabla
w_\epsilon|^2$ is bounded from above uniformly in $\epsilon$,
which implies that $w\in {\cal D}^{1,2}(\R^N)$. Then, we can use
the second part of Lemma~\ref{w-eps} and choose a $w_\epsilon$ so
that, additionally,
$\|w-w_\epsilon\|_{\mathcal{D}^{1,2}}\le\epsilon$.

Consequently,

\be \label{kappa_from_below0} \kappa_{m,n}\ge
\limsup\|v_k^\epsilon\|_{\eich_0}^2+(1-4\epsilon)\int|\nabla
w|^2.\ee

By assumption, $S_n>\kappa_{m,n}$, so there exists an $\epsilon>0$
such that $(1-4\epsilon)S_n>\kappa_{m,n}$. From
(\ref{kappa_from_below0}) then follows:

\be \label{kappa_from_below} \kappa_{m,n}>\kappa_{m,n}
t^\frac{2}{2^*}+\kappa_{m,n}(1-t)^\frac{2}{2^*},\ee for all
$t\in[0,1)$, which is false. Thus the assumption $t\neq 1$ is
false and  by (\ref{bll}), from $t=1$ follows $w=0$.
\end{proof}

\begin{remark}
\label{runaway1} Let $\Omega\subset\Domain$ be an open set. The
assertion of Lemma~\ref{runaway} holds also if $u_k$ is a
minimizing sequence for $\kappa_{m,n}(\Omega)$. The only
modification required for the proof is that inequality
$\kappa_{m,n}(\Omega)<S_N$ (due to Remark~\ref{kappa-S1}) replaces
$\kappa_{m,n}<S_N$.
\end{remark}

Let $\chi\in C_0^\infty((0,\infty))$ be the following even
function: $\chi(t)=0$ when $t\le\frac12$ or $t\ge 4$,
$\chi(t)=2(t-\frac12)$ when $t\in[\frac12,1]$, $\chi(t)=t$ when
$t\in[1,2]$, $\chi(t)=\frac12(4-t)$ when $t\in[2,4]$. Let
\begin{equation}
\chi_j=2^j\chi(2^{-j}|\cdot|), j\in\R.
\end{equation}

Let
\begin{equation}
\label{B} B_j=(0,2^{j})^{n}\times\{2^j<|y|<2^{j+1}\}, j\in\R.
\end{equation}

\begin{lemma}
\label{lem:step1} Let $u_k\in \eich_0$ be a bounded sequence. If
$u_k\stackrel{D_0}{\rightharpoonup}0$ and,for every sequence
$(x_k,y_k)\in\domain$, $u_k(\cdot+(x_k,y_k))\rightharpoonup 0$ in
$H^1_{loc}(\domain)$, then for every sequence $j_k\in\R$ \be
\label{step1} \int_{B_{j_k}}\chi_0(u_k)^{2^*}\to 0.\ee
\end{lemma}
\begin{proof}
It suffices to consider three cases: 1) $j_k\to-\infty$; 2)
$j_k\to +\infty$ and 3) $j_k$ is a bounded sequence.

Case 1. If $j_k\to-\infty$, $\int_{B_{j_k}}\chi_0(u_k)^{2^*}\le
C|{B_{j_k}}|\to 0$.

Case 2. Assume that $|j_k|\le M\in\R$. Then
$\int_{B_{j_k}}\chi_0(u_k)^{2^*}\to 0$ since $u_k\rightharpoonup
0$ in ${\cal D}^{1,2}_{loc}$, and in particular, $u_k\to 0$
locally in measure.

Case 3. Assume that $j_k\to +\infty$ and, without loss of
generality, that $j_k\in\N$. Consider a tesselation of $B_{j_k}$
by the sets $B_{j_k}^{il}=Q_i\times\{2^{j_k}+l<|y|<2^{j_k}+l+1\}$,
where $Q_i$ are unit cubes in $\R^n$ and $l=0,\dots,2^{j_k+1}-1$.
We will use the following version of the Sobolev inequality that
holds for all $i,l$ with a uniform constant $C$:

\be \label{uglySob}
\left(\int_{B_{j_k}^{il}}w^{2^*}\right)^\frac{2}{2^*}\le
C\left(\int_{B_{j_k}^{il}} |y|^{2-m}|\nabla(|y|^\frac{m-2}{2}w)|^2
+ \int_{B_{j_k}^{il}} |y|^{-2}|w|^2\right).\ee

Substituting $w=\chi_0(u_k)$ and taking into account that
$\chi_0(t)^2\le C t^{2^*}$ and that $|y|^{-2}\le 2^{-j}$, we have,
with a renamed constant,

\be \int_{B_{j_k}^{il}}\chi_0(u_k)^{2^*}\le
C\left(\int_{B_{j_k}^{il}}
|y|^{2-m}|\nabla(|y|^\frac{m-2}{2}u_k)|^2 + \int_{B_{j_k}^{il}}
|u_k|^{2^*}\right)
\left(\int_{B_{j_k}^{il}}\chi_0(u_k)^{2^*}\right)^{1-\frac{2}{2^*}}
.\ee

Adding the inequalities above over all $i,l$, we get

\be \int_{B_{j_k}}\chi_0(u_k)^{2^*}\le C\left(\int_{B_{j_k}}
|y|^{2-m}|\nabla(|y|^\frac{m-2}{2}u_k)|^2 + \int_{B_{j_k}}
|u_k|^{2^*}\right)\sup_{i,l}\left(\int_{B_{j_k}^{il}}\chi_0(u_k)^{2^*}\right)^{1-\frac{2}{2^*}}
.\ee

Note that the first factor in the right hand side is bounded by
$\|u_k\|_{\eich_0}$. Hence, in order to verify (\ref{step1}) it
suffices to show that

\be \label{lsl} \int_{B_0}\chi_0(u_k(\cdot-(x_k,y_k)))^{2^*}\to
0\ee

for all $(x_k,y_k)\in\domain$ with $|y_k|\to\infty$. Indeed, since
by assumption $u_k(\cdot-(x_k,y_k))\rightharpoonup 0$ in
$H^1(B_0)$,

\be \int_{B_0}|u_k(\cdot-(x_k,y_k))|^2\to 0 \mbox{ for
}|y_k|\to\infty.\ee
>From here follows (\ref{lsl}), and therefore,
(\ref{step1}), once we take into account that $\chi(t)^{2^*}\le
Ct^2$.
\end{proof}

\section{Existence of the minimizer}

We start this section in interpreting the conclusion of
Lemma~\ref{lem:step1} in terms of $\eich,D$. The subsequent proofs
will be carried out in the space $\eich$.

\begin{lemma}
\label{lem:step2} Let $v_k\in \eich$ be a bounded sequence such
that for all sequences $(x_k,y_k)\in\domain$, $t_k>0$, \be
t_k^\frac{N-2}{2}Tv_k(t_k\cdot+(x_k,y_k))\rightharpoonup 0.\ee
Then for every sequence $j_k\in\Z$ \be \label{step2}
\int_{B_{j_k}}|y|^\betastar\chi_0(v_k)^{2^*}\to 0.\ee
\end{lemma}
\begin{proof} Let $j\in\R$, $v\in\eich$, $u=Tv$. Then, observing
that $|y|\in (2^j,2^{j+1})$ on $B_j$,

\begin{eqnarray}
\label{cheq}
\int_{B_{j}}|y|^\betastar\chi_0(v)^{2^*}=\\
\nonumber\int_{B_{j}}|y|^\betastar\chi_0(|y|^\frac{m-2}{2}u)^{2^*}\le\\
\nonumber
C2^{j\betastar}\left(\int_{B_{j}}\chi_0(2^{j\frac{m-2}{2}}u)^{2^*}+
\int_{B_{j}}\chi_0(2^{(j+1)\frac{m-2}{2}}u)^{2^*}\right).
\end{eqnarray}

Let us estimate the first integral in the last expression. The
estimate of the second integral is totally analogous and may be
omitted. Let $t_j=2^{j\frac{m-2}{N-2}}$ and let
 $u^j=t_j^\frac{N-2}{2}u(t_j\cdot)$. Then

\begin{equation}
2^{j\betastar}\int_{B_{j}}\chi_0(2^{j\frac{m-2}{2}}u)^{2^*}\le\\
\int_{B_{\frac{nj}{N-2}}}\chi_0(u^j).
\end{equation}

Let now $j_k$ be an arbitrary sequence and substitute $j=j_k$,
$u=u_k:=Tv_k$:

\begin{equation}
\label{32}
2^{{j_k}\betastar}\int_{B_{j_k}}\chi_0(2^{{j_k}\frac{m-2}{2}}u_k)^{2^*}\le\\
\int_{B_{j_k\frac{n}{N-2}}}\chi_0(u_k^{j_k}),
\end{equation}
where $u_k^{j_k}:=t_{j_k}^\frac{N-2}{2}u(t_{j_k}\cdot)$ still
satisfies the assumptions of the lemma and thus the assumptions of
Lemma~\ref{lem:step1}. From the latter follows that the right hand
side in (\ref{32}) converges to zero, and tracing back
(\ref{cheq}) with $u=u_k=Tv_k$, $j=j_k$, we arrive at
(\ref{step2}).
\end{proof}

\begin{lemma}
\label{Lieb} If $v_k\in \eich$ is as in Lemma~\ref{lem:step2},
then $v_k\to 0$ in $L^{2^*}(\domain,\mymu)$.
\end{lemma}

\begin{proof} Let us use the following version of Sobolev inequality with a
fixed $q\in(2,2^*)$:

\be \label{uglySob2} C\left(\int_{B_0}|y|^\betastar
|v|^{2^*}\right)^\frac{q}{2^*}\le \int_{B_{0}} |y|^{2-m}|\nabla
v|^\frac{q}{2} + \int_{B_{0}}|y|^\betaq|w|^q.\ee The exponent
$\betaq$ is chosen so that the integral of the respective
expression over the whole $\Domain$ is dilation invariant, so that
the inequality holds (with the same constant) with $B_0$ replaced
with $\eta_{\alpha,j}B_0$, for all $j\in\Z$, $\alpha\in\R^n$.
Substituting $v=\chi_i(v_k)$, $i\in\Z$, we get

\begin{eqnarray} \label{uglySob3}
\int_{\eta_{\alpha,j}B_0}|y|^\betastar\chi_i(v_k)^{2^*}\le\\
\nonumber C\left( \int_{\eta_{\alpha,j}B_0\cap\supp \chi_i(v_k)}
|y|^{2-m}|\nabla v_k|^\frac{q}{2} +
\int_{\eta_{\alpha,j}B_0}|y|^{\betaq}\chi_i(v_k)^q\right)\times
\left(\int_{\eta_{\alpha,j}B_0}\chi_i(v_k)^{2^*}\right)^{1-\frac{q}{2^*}}.\end{eqnarray}

Adding terms up over all $i,j\in\Z$ and $\alpha\in\Z^n$, we have

\begin{eqnarray*} \label{uglySob4} \int_{\domain}|y|^\betastar|v_k|^{2^*}\le
C\left(\int_{\domain} |y|^{2-m}|\nabla v_k|^\frac{q}{2} +
\int_{\domain}|y|^{\betaq}|v_k|^q\right)\times\\\sup_{i,j,\alpha}
\left(\int_{\eta_{\alpha,j}B_0}|y|^\betastar\chi_i(v_k)^{2^*}\right)^{1-\frac{q}{2^*}}.
\end{eqnarray*}

Note that the first factor in the last expression is bounded,
since $v_k$ is a bounded sequence in $\eich$. In particular, \be
\int_{\domain}|y|^{\betaq}|v_k|^q \le C\|v_k\|_{\eich}^q\ee due to
the correspondent inequality (\cite{M}, p.98, Corollary 3. Thus
the lemma is proved once we verify that for an arbitrary sequence
$i_k,j_k\in\Z$ and $\alpha_k\in\Z^n$,

\be
\int_{\eta_{j_k,\alpha_k}B_0}|y|^\betastar\chi_{i_k}(v_k)^{2^*}\to
0. \ee This, however, is an immediate corollary of (\ref{step2}),
once we substitute $v_k=t_k^\frac{n}{2}{\tilde v}_k(t_k\cdot)$
with a suitable sequence $t_k$.
\end{proof}

\begin{corollary}
\label{cor:novanish} Let $v_k$ be a minimizing sequence for
(\ref{HSv}), namely, $\|v_k\|^2_\eich\to\kappa_{m,n}$,
$\|v_k\|_{2^*;\mymu}=1$. Then there is a sequence $g_k\in D$, such
that, on a renamed subsequence, $\wlim g_kv_k\neq 0$.
\end{corollary}

\begin{proof} Assume the opposite, namely that $v_k\cw 0$. Note
that $v_k$ is a minimizing sequence, and so is $g_kv_k$ with any
$g_k\in D$. Then by Lemma~\ref{runaway} the sequence $u_k=Tv_k$
satisfies the assumptions of Lemma~\ref{Lieb} and thus $v_k\to 0$
in $L^{2^*}(\domain,\mymu)$, a contradiction.
\end{proof}

We now can prove Theorem~\ref{main}, from which
Theorem~\ref{main2} follows immediately due to the isometry
(\ref{transform}).

{\em Proof of Theorem~\ref{main}.} Let $v_k$ be a minimizing
sequence. Due to Corollary~\ref{cor:novanish}, we may assume
without loss of generality that $v_k\rightharpoonup w\neq 0$. Then
$t:=\int_{\domain} |w|^{2^*}\mymu\in (0,1]$. From Brezis-Lieb
lemma follows then that $\int_{\domain} |v_k-w|^{2^*}\mymu=1-t$.
Therefore, \be
\kappa_{m,n}=\lim\|v_k\|^2_\eich=\|w\|^2_\eich+\lim\|v_k-w\|_\eich^2\ge
\kappa_{m,n}t^\frac{2}{2^*}+\kappa_{m,n}(1-t)^\frac{2}{2^*}.\ee
This inequality holds only as equality at the endpoints $t=0,1$
and thus, with necessity, $t=1$. In other words, $v_k\to w$ in
$L^{2^*}(\domain,\mymu)$ and therefore $w$ is a minimizer. \qed

\section{Existence and non-existence for $\Omega\subset\Domain$}

In this section we give the proofs of Theorem~\ref{main6} and
Theorem~\ref{main4}.

{\em Proof of Theorem~\ref{main6}.} Let $u_\epsilon\in
C_0^\infty(\Domain)$, $\epsilon>0$, satisfy
\be\label{approx}\|u\|_{2^*,|y|^\betastar}=1 \mbox{ and }
\|u_\epsilon\|^2_{\eich}\le\kappa_{m,n}+\epsilon \ee and set
$v_{\epsilon,t}:=t^\frac{N-2}{2}u_\epsilon(t\cdot)$, $t>0$.

Case (a): Without loss of generality assume that $x_0=0$.  Then
for $t>0$ sufficiently large $v_{\epsilon,t}\in
C_0^\infty(B_r(0))$ and still satisfies (\ref{approx}).
Consequently, $\kappa_{m,n}(B_r(0))\le\kappa_{m,n}+\epsilon$, and
since $\epsilon$ is arbitrary,
$\kappa_{m,n}(B_r(0)))\le\kappa_{m,n}$. The converse inequality
$\kappa_{m,n}(B_r(0))\ge\kappa_{m,n}$ is immediate.

Case (b): The proof is completely analogous to the case (a) once
we note that for $t>0$ sufficiently small $v_{\epsilon,t}\in
C_0^\infty(\{|y|>R\})$. \qed\vskip3mm

We proceed now with the proof of Theorem~\ref{main4}. The
following statement is a particular case of the global compactness
theorem from \cite{acc} to the case $\TheSpace(\R^N)$ with the
group $D_1$ of unitary operators generated by actions of dilations

\be (h_tu)(x)=t^\frac{N-2}{2}u(tx),\,t>0,\ee

and translations

\be u\mapsto u(\cdot-z), z\in\R^N. \ee

\begin{theorem}
\label{abstractcc}
  Let $u_{k}  \in \TheSpace(\R^N)$ be a  bounded sequence. Then there exists
$w^{(\ell)}  \in \TheSpace(\R^N)$, $g_{k} ^{(\ell)}  \in D_1$,
$k,\ell\in\N$, such that for a renumbered subsequence one has:
\begin{eqnarray}
\label{separates} g_k^{(1)}&=&id,\; {g_{k} ^{(i)}} ^{-1}  g_{k}
^{(j)}\rightharpoonup 0 \mbox{ for } i\neq j,
\\
w^{(\ell)}&=&\wlim  {g_{k} ^{(\ell)}}^{-1}u_k
\\
\label{norms} \sum_{\ell\in\N} \|w ^{(\ell)}\|^2 & \le & \limsup
\|u_k\|^2
\\
\label{BBasymptotics} u_{k}  - \sum_{\ell\in\N}  g_{k} ^{(\ell)}
w^{(\ell)} &\to& 0 \mbox{ in } L^{2^*}(\R^N).
\end{eqnarray}
The series in (\ref{BBasymptotics}) is absolutely convergent in
$\TheSpace$, uniformly in $k$.
\end{theorem}

\begin{lemma}
\label{BLsum} Let $u_k$, $w^{(i)}$ be as in
Theorem~\ref{abstractcc}. Then

\be \label{BLsum-Omega} \int_{\R^N} |u_k|^{2^*}\to\sum_i \int_
{\R^N}|w^{(i)}|^{2^*}.\ee
\end{lemma}

\begin{proof} By (\ref{BBasymptotics}) and continuity of
$u\mapsto \int_{\R^N} |u|^{2^*}$ in $\TheSpace(\R^N)$ it suffices
to prove the lemma for $u_k^M:=\sum_{i=1}^Mg_k^{(i)}w^{(i)}$,
$M\in\N$. Iterating Brezis-Lieb lemma for $M-1$ steps, we obtain
immediately

\be  \int_{\R^N} |u_k^M|^{2^*}=\sum_{i=1}^M \int_
{\R^N}|w^{(i)}|^{2^*}.\ee
\end{proof}

\begin{lemma}
\label{supports} Let $\Omega\subsetneq\R^N$ be an open set with
$\partial\Omega\in C^1$ and let $u_k\in H_0^1(\Omega)$. If there
exist $t_k>0$, $z_k\in\R^N$, $w\in \TheSpace(\R^N)$ such that
$t_k^{-\frac{N-2}{2}}u_k(t_k^{-1}\cdot+z_k)\rightharpoonup w$ in
$\TheSpace(\R^N)$, then, modulo a set of measure zero, \be V(w):=
\{z\in\R^N:w(z)\neq 0\}\subset\liminf t_k(\Omega-z_k).  \ee
Moreover, if there exist $t_0>0$ and $z_0\in\R^N$ such that
$\liminf t_k(\Omega-z_k)\subset t_0(\Omega+z_0)$, then
$w(t_0^{-1}(\cdot-z_0))\in H_{0,loc}^1(\Omega)$.
\end{lemma}

\begin{proof} Convergence $t_k^{-\frac{N-2}{2}}u_k(t_k^{-1}\cdot+z_k)\rightharpoonup w$
in $\TheSpace(\R^N)$ implies convergence a.e. Therefore, in the
complement of a set of measure zero, $w(z)=0$ for any $z$ that is
not in $\cap_{k\ge k_0}(t_k(\Omega-z_k))$ for some $k_0\in\N$.
Consequently, $w(z)=0$ unless, modulo a set of zero measure,
$z\in\liminf(t_k(\Omega-z_k))$. Then, by assumption,
$V(w(t_0^{-1}(\cdot-z_0))\subset\Omega$. Since $\partial\Omega\in
C^1$ and $w\in\TheSpace(\R^N)$, the conclusion
$w(t_0^{-1}(\cdot-z_0))\in H_{0,loc}^1(\Omega)$ follows.
\end{proof}

\noindent {\em Proof of Theorem~\ref{main4}}.

1. Observe that the norms $\eich_0(\Omega)$, $\TheSpace(\Omega)$
and $H^1_0(\Omega)$ are equivalent. The last two are equivalent by
the Friedrichs inequality, which holds since $\partial\Omega\in
C^1$ and $\sup_{x\in\Omega} d(x,\R^N\setminus\Omega)<\infty$.
Since $\|u\|_{\eich_0(\Omega)}\le\|u\|_{\TheSpace(\Omega)}$, it
suffices to show that

\be C(\Omega):=\inf_{u\in
C_0^\infty(\Omega):\;\int_\Omega\frac{u^2}{|y^2|}=1}\int_\Omega|\nabla
u|^2>\HardyC. \ee Let $R>0$ be such that $\Omega\subset
A_{1/R,R}$. Then it suffices to show that $C(A_{1/R,R})>\HardyC$.
This easily follows from \be \inf_{u\in
C_0^\infty(\omega_R):\;\int_{\omega_R}\frac{u^2}{|y^2|}dy=1}
\int_{\omega_R}|\nabla_y u|^2dy>\inf_{u\in
C_0^\infty(\R^m):\;\int_{\R^m}\frac{u^2}{|y^2|}dy=1}
\int_{\R^m}|\nabla_y u|^2dy=\HardyC,\ee where
$\omega_R=A_{1/R,R}\cap\R^m$, while the latter inequality holds
true since the minimum in the left hand side is attained (the
Dirichlet problem on a bounded domain) and the minimizer cannot be
a minimizer on $\R^m$ by the maximum principle.

2. Let $u_k$ be a minimizing sequence for (\ref{minProb4}).
By the preceding step $u_k$ is bounded in $\TheSpace(\R^N)$.
Assume that, on a renumbered subsequence,
\be
\label{wlw}
{t_k}^{-\frac{N-2}{2}}u(t_k\cdot+z_k)\rightharpoonup w\neq 0, \; z_k=(x_k,y_k), t_k>0
\ee
Note that if $t_k\to 0$, the
scaling argument gives ${t_k}^{-\frac{N-2}{2}}u_k(t_k\cdot+z_k)\to 0$ in
$L^2(\R^N)$, so $w=0$.
If, on the other hand, there is a subsequence where
both $t_k$ and $1/t_k$ are bounded, but $y_k$ is unbounded, from Lemma~\ref{supports} follows that the set $\{w\neq 0\}$
has measure zero, which also yields $w=0$.
We conclude that $w\neq 0$ only if
either (a) $t_k\to\infty$
or (b) $t_k$, $1/t_k$
and $y_k$ are bounded.

3. Let us show now that case (a) does not occur.
Assume that there is a sequence $(z_k,t_k)\in
\R^N\times(0,\infty)$, $t_k\to\infty$, such that, on a renamed subsequence,
$t_k^{-\frac{N-2}{2}}u_k(t_k^{-1}\cdot+z_k)\rightharpoonup w$ with
some $w\in{\mathcal D}^{1,2}(\R^N)\setminus\{0\}$.
%
Since $C_0^\infty(\R^N)$ is dense in ${\mathcal D}^{1,2}(\R^N)$,
for every $\epsilon\in(0,\frac12)$ there exists a $w_\epsilon\in
C_0^\infty(\R^N)$ such that

\begin{equation}
\label{weps} \|w_\epsilon-w\|_{{\mathcal
D}^{1,2}}\le\epsilon\|w\|_{{\mathcal D}^{1,2}}.
\end{equation}

Moreover, we can choose $w_\epsilon$ so that for all
$k$ sufficiently large $w_\epsilon(t_k(\cdot-z_k))$ is supported
in $\Omega+B_{2\epsilon}(0)$. Indeed, $w_\epsilon(t_k(\cdot-z_k))$
is supported in an arbitrarily small (for $k$ large) neighborhood
of $z_k$ and, since for every $z\in\R^N$, $u_k(t_k^{-1}z+z_k)=0$
whenever $z_k\notin\Omega+B_\epsilon(0)$ and $k$ is
sufficiently large, we have necessarily $z_k\in
\Omega+B_\epsilon(0)$.

Let
\begin{equation}
v_{\epsilon,k}=u_k-t_k^{\frac{N-2}{2}}w_\epsilon(t_k(\cdot-z_k)).
\end{equation}
Then, using the scaling invariance of the involved integrals and reserving the norm notation
for the $\eich_0$-norm, one has

\begin{eqnarray*}
\|u_k\|^2 =
\|v_{\epsilon,k}+t_k^{\frac{N-2}{2}}w_\epsilon(t_k(\cdot-z_k))\|^2&=&
\\
\|v_{\epsilon,k}\|^2+
\|t_k^{\frac{N-2}{2}}w_\epsilon(t_k(\cdot-z_k))\|^2+
2\int\nabla v_{\epsilon,k}\cdot\nabla
t_k^{\frac{N-2}{2}}w_\epsilon(t_k(\cdot-z_k)))&-&
\\
2\HardyC \int
v_{\epsilon,k}t_k^{\frac{N-2}{2}}w_\epsilon(t_k(\cdot-z_k)))|y|^{-2}&=&
\\
\|v_{\epsilon,k}\|^2+\int|\nabla w_\epsilon|^2+o(1) +
2\int\nabla u_k\cdot\nabla
t_k^{\frac{N-2}{2}}w_\epsilon(t_k(\cdot-z_k))-2\int|\nabla
w_\epsilon|^2& -&
\\
2\HardyC\int
u_kt_k^{\frac{N-2}{2}}w_\epsilon(t_k(\cdot-z_k))|y|^{-2}+ 2\HardyC\int
t_k^{N-2}|w_\epsilon(t_k(\cdot-z_k))|^2|y|^{-2}&=&
\\
\|v_{\epsilon,k}\|^2+\int|\nabla w_\epsilon|^2+
2\int\nabla
w\cdot\nabla w_\epsilon-2\int|\nabla w_\epsilon|^2+o(1).
\end{eqnarray*}
At the last step we have used the following estimates:

\begin{equation}
\int t_k^{N-2}|w_\epsilon(t_k(\cdot-z_k))|^2|y|^{-2}=
t_k^{-2}\int|w_\epsilon|^2|y|^{-2}=o(1)
\end{equation}
and, by Cauchy inequality,

\begin{equation}
\left|\int
u_kt_k^{\frac{N-2}{2}}w_\epsilon(t_k(\cdot-z_k))|y|^{-2}\right|\le
\\
\left(\int u_k^2|y|^{-2}\right)^\frac12
\left(\int t_k^{N-2}|w_\epsilon(t_k(\cdot-z_k))|^2|y|^{-2}\right)^\frac12=o(1).
\end{equation}
Consequently, using (\ref{weps}), we obtain

\begin{equation}
\kappa_{m,n}(\Omega)=\lim \|u_k\|^2\ge \liminf
\|v_{\epsilon,k}\|^2+(1-8\epsilon)\int|\nabla w|^2.
\end{equation}
From the definitions of $\kappa_{m,n}(\Omega)$ and $S_N$ then follows

\begin{equation}
\label{kl} \kappa_{m,n}(\Omega)\ge \kappa_{m,n}(\Omega)\liminf
\|v_{\epsilon,k}\|_{2^*}^2+(1-8\epsilon)S_N\|w\|_{2^*}^2.
\end{equation}

Let now
\begin{equation}
v_{k}=u_k-t_k^{\frac{N-2}{2}}w(t_k(\cdot-z_k)).
\end{equation}
Then (\ref{weps}) and (\ref{kl}) imply

\begin{equation}
 \kappa_{m,n}(\Omega)\ge \kappa_{m,n}(\Omega)\liminf
\|v_{k}\|_{2^*}^2+(1-10\epsilon)S_N\|w\|_{2^*}^2,
\end{equation}
and, since $\epsilon$ is arbitrary,

\begin{equation}
\label{BNenergy} 1\ge \liminf
\|v_{k}\|_{2^*}^2+\frac{S_N}{\kappa_{m,n}(\Omega)}\|w\|_{2^*}^2.
\end{equation}
At the same time, from the Brezis-Lieb lemma \cite{BL} (passing to a
renamed subsequence if necessary) follows

\begin{equation}
\label{BNmass} \liminf \|v_{k}\|_{2^*}^{2^*}+\|w\|_{2^*}^{2^*}=1.
\end{equation}
Since $S_N>\kappa_{m,n}(\Omega)$ by Remark~\ref{kappa-S1} and $2^*>2$,
relations (\ref{BNmass}) and (\ref{BNenergy}) hold simultaneously only
if $w=0$, a contradiction.

4. We conclude that case (b) is the only possibility for a non-zero weak limit
(\ref{wlw}). Consequently, for a renumbered subsequence, one can write (\ref{BBasymptotics})  as
\be
u_k-\sum w^{(i)}(\cdot-x_k^{(i)})\to 0 \text { in } L^p(\R^N), p\in(2,2^*],
\ee
$|x_k^{(i)}-x_k^{(j)}|\to\infty$ for $i\neq j$,
and, moreover (by using (\ref{flask-x}), Lemma~\ref{supports} and the fact that $\TheSpace(\Omega))$-norm and the $\eich_0(\Omega)$) are equivalent,
$w^{(j)}\in H^1_0(\Omega)$.
>From Lemma~\ref{BLsum} follows that
\be
\label{0mass}
\int|u_k|^{2^*}\to \sum_j\int|w^{(j)}|^{2^*}=1,
\ee
while (\ref{norms}) and definition (\ref{minProb4a}) of $\kappa_{m,n}(\Omega)$ imply

\be
\label{0nrg}
\kappa_{m,n}(\Omega)\ge \sum_j \kappa_{m,n}(\Omega)\left(\int|w^{(j)}|^{2^*}\right)^{2/2^*}.
\ee
Relations (\ref{0mass}) and (\ref{0nrg}) are contradictive unless all but one $w^{(j)}$ equal zero and for some $j_0$, $\int|w^{(j_0)}|^{2^*}=1$. Then, necessarily, $w^{(j_0)}$ is a minimizer.
\qed

\section{Existence in $R^{n+m}$ - concluding remarks}

One can use the rearrangement argument to reduce the proof of Theorem~\ref{main} in the case $m>2$ to an existence result of Badiale and Tarantello  \cite{BadialeTarantello}. It should be noted that this reduction does not extend to the case $m=1$, while the proof in the present paper is uniform with regard to $m$ and uses techniques that allow to approach analogous problems with lack of radial symmetry.

\begin{proof}
By the rearrangement argument, the minimum in (\ref{minProb})
is attained in the subspace $\eich_r$ of $\eich$ of functions that are radially symmetric in the variable
$y\in\R^n$, so we may restate the problem, regarding $|y|$ as a
radial variable in $\R^2$, in the form

\begin{equation}
\label{CNNr} \kappa_{m,n}=\inf_{v\in \mathcal{D}_r^{1,2}(\R^{n+2}):
\frac{\omega_{m-1}}{2\pi}\int_{R^n\times\R^2}|v(x,y)|^\frac{2N}
{N-2}|y|^{-\frac{2(m-2)}{N-2}}dxdy=1}
\frac{\omega_{m-1}}{2\pi}\int_{\R^n\times\R^2}|\nabla
u(x,y)|^2dxdy.
\end{equation}

We note that the exponent $\frac{2N}{N-2}$ is subcritical in the dimension $n+2<N$
and we may apply Theorem~2.5 of \cite{BadialeTarantello}, relative to $\R^{n+2}$. Parameters for the application,
in the original notations of \cite{BadialeTarantello}, are $q=2$,
$s=\frac{2(m-2)}{N-2}$ and $q_*=\frac{2N}{N-2}$.
\end{proof}

Several related problems remain unresolved in this paper.

{\bf 1.} Evaluate the best constant $\kappa_{m,n}$ and find
minimizers of (\ref{HS}) when the minimum exists. Is the
inequality $\kappa_{1,2}< S_{3}$ true? If it is false, is there
still a minimizer for $\kappa_{1,2}$?

{\bf 2.} We saw that when $\kappa_{m,n}(\Omega)=\kappa_{m,n}$,
$\Omega\neq\Domain$, there is no minimizer for
$\kappa_{m,n}(\Omega)$. Is the converse true? That is, does
$\kappa_{m,n}(\Omega)>\kappa_{m,n}$ imply existence of the
minimizer?

{\bf 3.} More general results of the form (\ref{HS}) were recently
established (\cite{FMT},\cite{FMT1}). For example if $\xO$ is a
bounded smooth and convex domain, $d(x) = {\rm dist}(x, \partial
\xO)$ then there exist a   positive constant $C$ dependent on $\xO$ such that
\be\la{25} \int_{\xO} |\nabla u|^2 dx-  \frac14 \int_{\xO}
\frac{u^2}{d^2} dx
  \geq
 C \left( \int_{\xO} |u|^{\frac{2N}{N-2}} dx \right)^{\frac{N-2}{N}},
\qquad \quad  \forall u \in C_0^{\infty}(\xO). \ee We believe that
the following is an interesting question. Is the best
 constant $C=C(\xO)$
connected with the constant $\kappa_{1,n}$? In particular, is it
true for convex $\Omega$ that $C(\Omega)=\kappa_{1,n}$?

{\bf 4.} By analogy with Theorem~\ref{main}, it is natural to ask
whether the following minimization problem
\begin{equation}
\label{minProb3}\kappa_{m,n}^p=\inf_{\int_{\domain}|y|^{\frac{N(p-m)}{N-p}}|v|^{\frac{Np}{N-p}}dxdy=1}\int_{\Domain}
|y|^{-(m-p)}|\nabla v|^pdxdy
\end{equation}
has a minimizer when $m\neq p>1$, in particular when $p^2<N$.
A partial answer can be given by an argument similar to the previous section, by Theorem~2.5 of \cite{BadialeTarantello}.

\end{document}